# Recent Developments in Nonparametric Inference and Probability

## Festschrift for Michael Woodroofe

Jiayang Sun, Anirban DasGupta, Vince Melfi, Connie Page,
Editors

# Contents

# Preface

There have been extensive developments recently in modern nonparametric inference and modeling. Nonparametric and semi-parametric methods are especially useful with large amounts of data that are now routinely collected in many areas of science. Probability and stochastic modeling are also playing major new roles in scientific applications.

This refereed special volume highlights challenges and developments at this interface of statistics, probability and the sciences, and honors Michael B. Woodroofe for his pioneering contributions to nonparametric inference and probability.

Topics covered include biased sampling and missing data, shape-restricted inference, contemporary sequential analysis, modern nonparametric inference, probability, and statistics applications. Some of the papers in this volume were presented at "A Conference on Nonparametric Inference and Probability with Applications to Science", September 24–25, 2005, at Ann Arbor, Michigan:

http:// www.stat.lsa.umich.edu/conference/mw2005/

for which there were 97 pre-registered participants. The conference was organized by Robert Keener and Jiayang Sun. Invited Speakers and participants included: Persi Diaconis, Robert Keener, Steve Lalley, T.L. Lai, Mario Mateo, Mary Meyer, Vijay Nair, Gordie Simons, Jiayang Sun, Michael Woodroofe, Wei Biao Wu, Cun-Hui Zhang, Charles Hagwood, Steve Coad, Hira Koul, Anand Vidyashenkar, Zhiliang Ying, Arniban DasGupta, Moulinath Banerjee, Anna Amirdjanova, Bill Strawderman, Connie Page, Tom Sellke, Byron Roe, Vince Melfi, and Herman Chernoff.

Michael Woodroofe has had a distinguished career and is widely recognized as a pre-eminent scientist in statistics and probability. He has broad interests and has made deep and significant contributions in many areas, many represented by sessions at this conference. His work in probability, sequential analysis and inference is regarded as seminal and pioneering and has influenced an entire generation of researchers. He has published more than 100 research articles, written a SIAM monograph, and authored a book. He is a former Editor of the Annals of Statistics, a member of Phi Beta Kappa, and a fellow of the Institute of Mathematical Statistics and an elected member of the International Statistical Institute.

We thank the National Science Foundation, the National Security Agency, and the University of Michigan, (the Department of Statistics; the College of Literature, Science, and the Arts; the Office of the Vice President for Research; and the Rackham Graduate School) for their financial support to the conference, and the Institute of Mathematical Statistics, and the American Statistical Association for the co-sponsorship. We also thank Amy Rundquist for conference administration.

Jiayang Sun,
*Case Western Reserve University*
Vince Melfi,
*Michigan State University*

Anirban DasGupta,
*Purdue University*
Connie Page,
*Michigan State University*

## Special Thanks to Editorial Assistants

Steve Ganocy
Xiaofeng Wang

# Contributors to this volume

Bartroff, J., *Stanford University*
Bogie, K., *Cleveland FES center*

Chuang, C. S., *Millennium Partners*
Coad, D. S., *University of London*
Cohen, A., *Rutgers University*

DasGupta, A., *Purdue University*

Heller, B., *Illinois Institute of Technology*

Keener, R., *University of Michigan*
Kim, D.-Y., *Michigan State University*
Kolassa, J., *Rutgers University*

Lai, T. L., *Stanford University*

Marchand, E., *Université de Sherbrooke*

Nagai, K., *Yokohama National University*

Sackrowitz, H., *Rutgers University*
Sarkar, J., *Indiana University Purdue University Indianapolis*
Shao, X., *University of Chicago*
Strawderman, W. E., *Rutgers University*
Su, Z., *Stanford University*
Sun, J., *Case Western Reserve University*

Wang, B., *University of South Alabama*
Wang, M., *University of Chicago*
Wang, X., *The Cleveland Clinic Foundation*
Weng, R. C., *National Chengchi University*
Woodroofe, M., *the University of Michigan*
Wu, W. B., *University of Illinois, at Urbana-Champaign*

Zhang, C.-H., *Rutgers University*
Zhang, T., *Purdue University*